\documentstyle[diagram,12pt]{article}

\expandafter\chardef\csname pre amssym.def at\endcsname=\the\catcode`\@
\catcode`\@=11

\def\undefine#1{\let#1\undefined}
\def\newsymbol#1#2#3#4#5{\let\next@\relax
 \ifnum#2=\@ne\let\next@\msafam@\else
 \ifnum#2=\tw@\let\next@\msbfam@\fi\fi
 \mathchardef#1="#3\next@#4#5}
\def\mathhexbox@#1#2#3{\relax
 \ifmmode\mathpalette{}{\m@th\mathchar"#1#2#3}%
 \else\leavevmode\hbox{$\m@th\mathchar"#1#2#3$}\fi}
\def\hexnumber@#1{\ifcase#1 0\or 1\or 2\or 3\or 4\or 5\or 6\or 7\or 8\or
 9\or A\or B\or C\or D\or E\or F\fi}

\font\tenmsa=msam10 scaled \magstep1
\font\sevenmsa=msam7 scaled \magstep1
\font\fivemsa=msam5 scaled \magstep1
\newfam\msafam
\textfont\msafam=\tenmsa
\scriptfont\msafam=\sevenmsa
\scriptscriptfont\msafam=\fivemsa
\edef\msafam@{\hexnumber@\msafam}
\mathchardef\dabar@"0\msafam@39
\def\dashrightarrow{\mathrel{\dabar@\dabar@\mathchar"0\msafam@4B}}
\def\dashleftarrow{\mathrel{\mathchar"0\msafam@4C\dabar@\dabar@}}

\def\ulcorner{\delimiter"4\msafam@70\msafam@70 }
\def\urcorner{\delimiter"5\msafam@71\msafam@71 }
\def\llcorner{\delimiter"4\msafam@78\msafam@78 }
\def\lrcorner{\delimiter"5\msafam@79\msafam@79 }
\def\yen{{\mathhexbox@\msafam@55 }}
\def\checkmark{{\mathhexbox@\msafam@58 }}
\def\circledR{{\mathhexbox@\msafam@72 }}
\def\maltese{{\mathhexbox@\msafam@7A }}

\font\tenmsb=msbm10 scaled \magstep1
\font\sevenmsb=msbm7 scaled \magstep1
\font\fivemsb=msbm5 scaled \magstep1
\newfam\msbfam
\textfont\msbfam=\tenmsb
\scriptfont\msbfam=\sevenmsb
\scriptscriptfont\msbfam=\fivemsb
\edef\msbfam@{\hexnumber@\msbfam}
\def\Bbb#1{\fam\msbfam\relax#1}
\def\widehat#1{\setboxz@h{$\m@th#1$}%
 \ifdim\wdz@>\tw@ em\mathaccent"0\msbfam@5B{#1}%
 \else\mathaccent"0362{#1}\fi}
\def\widetilde#1{\setboxz@h{$\m@th#1$}%
 \ifdim\wdz@>\tw@ em\mathaccent"0\msbfam@5D{#1}%
 \else\mathaccent"0365{#1}\fi}
\font\teneufm=eufm10
\font\seveneufm=eufm7
\font\fiveeufm=eufm5
\newfam\eufmfam
\textfont\eufmfam=\teneufm
\scriptfont\eufmfam=\seveneufm
\scriptscriptfont\eufmfam=\fiveeufm

\catcode`\@=\csname pre amssym.def at\endcsname

\expandafter\ifx\csname pre amssym.tex at\endcsname\relax \else \endinput\fi
\expandafter\chardef\csname pre amssym.tex at\endcsname=\the\catcode`\@
\catcode`\@=11
\newsymbol\boxdot 1200
\newsymbol\boxplus 1201
\newsymbol\boxtimes 1202
\newsymbol\square 1003
\newsymbol\blacksquare 1004
\newsymbol\centerdot 1205
\newsymbol\lozenge 1006
\newsymbol\blacklozenge 1007
\newsymbol\circlearrowright 1308
\newsymbol\circlearrowleft 1309
\undefine\rightleftharpoons
\newsymbol\rightleftharpoons 130A
\newsymbol\leftrightharpoons 130B
\newsymbol\boxminus 120C
\newsymbol\Vdash 130D
\newsymbol\Vvdash 130E
\newsymbol\vDash 130F
\newsymbol\twoheadrightarrow 1310
\newsymbol\twoheadleftarrow 1311
\newsymbol\leftleftarrows 1312
\newsymbol\rightrightarrows 1313
\newsymbol\upuparrows 1314
\newsymbol\downdownarrows 1315
\newsymbol\upharpoonright 1316
 \let\restriction\upharpoonright
\newsymbol\downharpoonright 1317
\newsymbol\upharpoonleft 1318
\newsymbol\downharpoonleft 1319
\newsymbol\rightarrowtail 131A
\newsymbol\leftarrowtail 131B
\newsymbol\leftrightarrows 131C
\newsymbol\rightleftarrows 131D
\newsymbol\Lsh 131E
\newsymbol\Rsh 131F
\newsymbol\rightsquigarrow 1320
\newsymbol\leftrightsquigarrow 1321
\newsymbol\looparrowleft 1322
\newsymbol\looparrowright 1323
\newsymbol\circeq 1324
\newsymbol\succsim 1325
\newsymbol\gtrsim 1326
\newsymbol\gtrapprox 1327
\newsymbol\multimap 1328
\newsymbol\therefore 1329
\newsymbol\because 132A
\newsymbol\doteqdot 132B
 
\newsymbol\triangleq 132C
\newsymbol\precsim 132D
\newsymbol\lesssim 132E
\newsymbol\lessapprox 132F
\newsymbol\eqslantless 1330
\newsymbol\eqslantgtr 1331
\newsymbol\curlyeqprec 1332
\newsymbol\curlyeqsucc 1333
\newsymbol\preccurlyeq 1334
\newsymbol\leqq 1335
\newsymbol\leqslant 1336
\newsymbol\lessgtr 1337
\newsymbol\backprime 1038
\newsymbol\risingdotseq 133A
\newsymbol\fallingdotseq 133B
\newsymbol\succcurlyeq 133C
\newsymbol\geqq 133D
\newsymbol\geqslant 133E
\newsymbol\gtrless 133F
\newsymbol\sqsubset 1340
\newsymbol\sqsupset 1341
\newsymbol\vartriangleright 1342
\newsymbol\vartriangleleft 1343
\newsymbol\trianglerighteq 1344
\newsymbol\trianglelefteq 1345
\newsymbol\bigstar 1046
\newsymbol\between 1347
\newsymbol\blacktriangledown 1048
\newsymbol\blacktriangleright 1349
\newsymbol\blacktriangleleft 134A
\newsymbol\vartriangle 134D
\newsymbol\blacktriangle 104E
\newsymbol\triangledown 104F
\newsymbol\eqcirc 1350
\newsymbol\lesseqgtr 1351
\newsymbol\gtreqless 1352
\newsymbol\lesseqqgtr 1353
\newsymbol\gtreqqless 1354
\newsymbol\Rrightarrow 1356
\newsymbol\Lleftarrow 1357
\newsymbol\veebar 1259
\newsymbol\barwedge 125A
\newsymbol\doublebarwedge 125B
\undefine\angle
\newsymbol\angle 105C
\newsymbol\measuredangle 105D
\newsymbol\sphericalangle 105E
\newsymbol\varpropto 135F
\newsymbol\smallsmile 1360
\newsymbol\smallfrown 1361
\newsymbol\Subset 1362
\newsymbol\Supset 1363
\newsymbol\Cup 1264
 
\newsymbol\Cap 1265
 
\newsymbol\curlywedge 1266
\newsymbol\curlyvee 1267
\newsymbol\leftthreetimes 1268
\newsymbol\rightthreetimes 1269
\newsymbol\subseteqq 136A
\newsymbol\supseteqq 136B
\newsymbol\bumpeq 136C
\newsymbol\Bumpeq 136D
\newsymbol\lll 136E
 
\newsymbol\ggg 136F
 
\newsymbol\circledS 1073
\newsymbol\pitchfork 1374
\newsymbol\dotplus 1275
\newsymbol\backsim 1376
\newsymbol\backsimeq 1377
\newsymbol\complement 107B
\newsymbol\intercal 127C
\newsymbol\circledcirc 127D
\newsymbol\circledast 127E
\newsymbol\circleddash 127F
\newsymbol\lvertneqq 2300
\newsymbol\gvertneqq 2301
\newsymbol\nleq 2302
\newsymbol\ngeq 2303
\newsymbol\nless 2304
\newsymbol\ngtr 2305
\newsymbol\nprec 2306
\newsymbol\nsucc 2307
\newsymbol\lneqq 2308
\newsymbol\gneqq 2309
\newsymbol\nleqslant 230A
\newsymbol\ngeqslant 230B
\newsymbol\lneq 230C
\newsymbol\gneq 230D
\newsymbol\npreceq 230E
\newsymbol\nsucceq 230F
\newsymbol\precnsim 2310
\newsymbol\succnsim 2311
\newsymbol\lnsim 2312
\newsymbol\gnsim 2313
\newsymbol\nleqq 2314
\newsymbol\ngeqq 2315
\newsymbol\precneqq 2316
\newsymbol\succneqq 2317
\newsymbol\precnapprox 2318
\newsymbol\succnapprox 2319
\newsymbol\lnapprox 231A
\newsymbol\gnapprox 231B
\newsymbol\nsim 231C
\newsymbol\ncong 231D
\newsymbol\diagup 231E
\newsymbol\diagdown 231F
\newsymbol\varsubsetneq 2320
\newsymbol\varsupsetneq 2321
\newsymbol\nsubseteqq 2322
\newsymbol\nsupseteqq 2323
\newsymbol\subsetneqq 2324
\newsymbol\supsetneqq 2325
\newsymbol\varsubsetneqq 2326
\newsymbol\varsupsetneqq 2327
\newsymbol\subsetneq 2328
\newsymbol\supsetneq 2329
\newsymbol\nsubseteq 232A
\newsymbol\nsupseteq 232B
\newsymbol\nparallel 232C
\newsymbol\nmid 232D
\newsymbol\nshortmid 232E
\newsymbol\nshortparallel 232F
\newsymbol\nvdash 2330
\newsymbol\nVdash 2331
\newsymbol\nvDash 2332
\newsymbol\nVDash 2333
\newsymbol\ntrianglerighteq 2334
\newsymbol\ntrianglelefteq 2335
\newsymbol\ntriangleleft 2336
\newsymbol\ntriangleright 2337
\newsymbol\nleftarrow 2338
\newsymbol\nrightarrow 2339
\newsymbol\nLeftarrow 233A
\newsymbol\nRightarrow 233B
\newsymbol\nLeftrightarrow 233C
\newsymbol\nleftrightarrow 233D
\newsymbol\divideontimes 223E
\newsymbol\varnothing 203F
\newsymbol\nexists 2040
\newsymbol\Finv 2060
\newsymbol\Game 2061
\newsymbol\mho 2066
\newsymbol\eth 2067
\newsymbol\eqsim 2368
\newsymbol\beth 2069
\newsymbol\gimel 206A
\newsymbol\daleth 206B
\newsymbol\lessdot 236C
\newsymbol\gtrdot 236D
\newsymbol\ltimes 226E
\newsymbol\rtimes 226F
\newsymbol\shortmid 2370
\newsymbol\shortparallel 2371
\newsymbol\smallsetminus 2272
\newsymbol\thicksim 2373
\newsymbol\thickapprox 2374
\newsymbol\approxeq 2375
\newsymbol\succapprox 2376
\newsymbol\precapprox 2377
\newsymbol\curvearrowleft 2378
\newsymbol\curvearrowright 2379
\newsymbol\digamma 207A
\newsymbol\varkappa 207B
\newsymbol\Bbbk 207C
\newsymbol\hslash 207D
\undefine\hbar
\newsymbol\hbar 207E
\newsymbol\backepsilon 237F
\catcode`\@=\csname pre amssym.tex at\endcsname

\newcommand{\ga}{\alpha}     
\newcommand{\gb}{\beta}      
   
\newcommand{\gd}{\delta}     
  
\newcommand{\gz}{\zeta}      
\newcommand{\gee}{\eta}      
\newcommand{\gth}{\theta}    
      
\newcommand{\gk}{\kappa}  
\newcommand{\gl}{\lambda}    
\newcommand{\gm}{\mu}

\newcommand{\gp}{\pi}        
\newcommand{\gr}{\rho}       
\newcommand{\gs}{\sigma}

\newcommand{\go}{\omega}

\newcommand{\setof}[2]{{\{\; #1 \; \vert \; #2 \; \} } }
\newcommand{\seq}[1]{{\langle #1 \rangle} }
\newcommand{\card}[1]{{\vert #1 \vert} }

\renewcommand{\models}{\vDash}
\newcommand{\powerset}{{\cal P}}
\newcommand{\dom}{{\rm dom}}

\newcommand{\crit}{{\rm crit}}

\newcommand{\cf}{{\rm cf}}

\newcommand{\lra}{\longrightarrow}

\newtheorem{definition}{Definition}
\newtheorem{theorem}{Theorem}
\newenvironment{proof}{\noindent{\bf Proof:}}{\nopagebreak\mbox{}\newline
 \makebox[\textwidth]{\hfill$\blacklozenge$}\par\bigskip}

\newcommand{\implies}{\Longrightarrow}

\catcode`\@=11
\def\@begintheorem#1#2{\rm \trivlist \item[\hskip \labelsep{\bf #1\ #2:}]}
\def\@opargbegintheorem#1#2#3{\rm \trivlist
      \item[\hskip \labelsep{\bf #1\ #2\ (#3):}]}
\catcode`\@=12

\newtheorem{lemma}{Lemma}

\newtheorem{corollary}{Corollary}

\newtheorem{claim}{Claim}

\newtheorem{fact}{Fact}

\newcommand{\commtriangle}[6]
{
\medskip
\[
\setlength{\dgARROWLENGTH}{6.0em}
\begin{diagram}
\node{#1} \arrow[2]{e,t}{#6} \arrow{se,b}{#4}  \node[2]{#3} \\
\node[2]{#2} \arrow{ne,r}{#5}
\end{diagram}
\]
\medskip
}
\newcommand{\FP}{{\Bbb P}}
\newcommand{\FQ}{{\Bbb Q}}
\newcommand{\FR}{{\Bbb R}}
\newcommand{\FS}{{\Bbb S}}

\newcommand{\F}{{\cal F}}

\title{A model in which every Boolean algebra has many subalgebras}

\author{James Cummings\\
        Hebrew University of Jerusalem
        \thanks{Research supported by a Postdoctoral Fellowship at the
        Mathematics Institute, Hebrew University of Jerusalem}\\
        \\
        Saharon Shelah\\
        Hebrew University of Jerusalem
        \thanks{Research partially supported by the Basic Research Fund of
        the Israel Academy of Science. Paper number 530.}}

\begin{document}

\baselineskip=16pt
\binoppenalty=10000
\relpenalty=10000

\maketitle

\begin{abstract}
   We show that it is consistent with ZFC (relative to large cardinals)
 that every infinite Boolean algebra $B$ has an irredundant subset $A$ such that
 $2^\card{A} = 2^\card{B}$. This implies in particular that $B$ has
 $2^\card{B}$ subalgebras. We also discuss some more general problems
 about subalgebras and free subsets of an algebra.

   The result on the number of subalgebras in a Boolean algebra
 solves a question of Monk from \cite{Monk}. The paper is intended to
 be accessible as far as possible to a general audience, in particular we
 have confined the more technical material  to a ``black box'' at the
 end. The proof involves a variation on Foreman and Woodin's model 
 in which GCH fails everywhere.
\end{abstract}

\section{Definitions and facts}

  In this section we give some basic definitions, and prove a couple of
 useful facts about algebras and free subsets. We refer the reader to
 \cite{Jech} for the definitions of set-theoretic terms. Throughout
 this paper we are working in ZFC set theory.

\begin{definition} Let $M$ be a set. $F$ is a {\bf finitary function on $M$\/}
 if and only if $F: M^n \lra M$ for some finite $n$.
\end{definition}

\begin{definition} ${\cal M}$ is an {\bf algebra\/} if and only if 
 ${\cal M} = (M, \cal F)$ where $M$ is a set, and $\cal F$ is a set
 of finitary functions on $M$ which is closed under composition and
 contains the identity function.

In this case we say ${\cal M}$ is an {\bf algebra on $M$.}
\end{definition}

\begin{definition} Let ${\cal M} = (M, {\cal F})$ be an algebra and let $N \subseteq M$.
Then $N$ is a {\bf subalgebra of $\cal M$\/} if and only if $N$ is closed under all
 the functions in $\cal F$. In this case there is a natural algebra structure
 on $N$ given by ${\cal N} = (N, {\cal F} \restriction {\cal N})$.
\end{definition}

\begin{definition} Let ${\cal M} = (M, {\cal F})$ be an algebra. Then
\[
    Sub({\cal M}) = \setof{N}{\hbox{$N$ is a subalgebra of ${\cal M}$}}.
\]
\end{definition}

\begin{definition}   Let ${\cal M} = (M, {\cal F})$ be an algebra. If $A \subseteq M$
then
\[
   Cl_{\cal M}(A) = \setof{F(\vec a)}{F \in {\cal F}, \vec a \in A}.
\]
 $Cl_{\cal M}(A)$ is the least subalgebra of ${\cal M}$ containing $A$.
\end{definition}

\begin{definition} Let ${\cal M} = (M, {\cal F})$ be an algebra. If $A \subseteq M$
 then $A$ is {\bf free\/} if and only if 
\[
    \forall a \in A \; a \notin Cl_{\cal M}(A - \{a\}).
\]
\end{definition}

\medskip

 In the context of Boolean algebras, free subsets are more usually referred
 to as {\bf irredundant.}
   We can use free sets to generate large numbers of subalgebras using the
 following well-known fact.

\medskip

\begin{fact}
\label{freefact}
 If $A \subseteq M$ is free for ${\cal M}$ then $\card{Sub(\cal M)} \ge 2^\card{A}$.
\end{fact}

\begin{proof} For each $B \subseteq A$ let $N_B = Cl_{\cal M}(B)$. It follows from the
 freeness of $A$ that $N_B \cap A = B$, so that the map
 $B \longmapsto N_B$ is an injection from $\powerset A$ into $Sub(\cal M)$.
\end{proof}

\medskip

 The next fact (also well-known) gives us one way of generating irredundant subsets in a Boolean
 algebra.

\medskip

\begin{fact} \label{BAfact} If $\gk$ is a strong limit cardinal and $B$ is a Boolean algebra
 with $\card{B} \ge \gk$ then $B$ has an irredundant subset of cardinality $\gk$.
\end{fact}

\begin{proof}  We can build by induction a sequence
 $\seq{b_\ga: \ga < \gk}$ such that $b_\gb$ is not above any element in
 the subalgebra generated by $\vec b \restriction \gb$. The point is
 that a subalgebra of size less than $\gk$ cannot be dense, since $\gk$
 is strong limit. It follows from the property we have arranged that
 we have enumerated an irredundant subset of cardinality $\gk$.
\end{proof}  

\medskip 

  We make the remark that the last fact works even if $\gk = \go$.
The next fact is a technical assertion which we will use when we
 discuss free subsets and subalgebras in the general setting.

\medskip

\begin{fact}
\label{singularfact}
 Let $\gm$ be a singular strong limit cardinal. Let ${\cal M} = (M, {\cal F})$ be an algebra
 and suppose that $\gm \le \card{M} < 2^\gm$ and $\card{\cal F} < \gm$. Then $Sub(\cal M)$
 has cardinality at least $2^\gm$.
\end{fact}

\begin{proof} Let $\gl = \cf(\gm) + \card{\cal F}$, then $2^\gl < \gm$ since $\gm$ is singular
 and  strong limit.
 Define
\[
     P = \setof{A \subseteq M}{\card{A} = \cf(\gm)}.
\]
 From the assumptions on $M$ and $\gm$,
\[
   2^\gm = \gm^{\cf{\gm}} \le \card{M}^{\cf{\gm}} = \card{P} \le 2^{\gm.\cf{\gm}} = 2^\gm,
\]
so that $\card{P} = 2^\gm$. 

   Observe that if $A \in P$ then $\card{ Cl_{\cal M}(A)} \le \gl$. So if we
 define an equivalence relation on $P$ by setting
\[
    A \equiv B \iff Cl_{\cal M}(A) = Cl_{\cal M}(B),
\]
 then the classes each have size at worst $2^\gl$. Hence there are
$2^\gm$ classes and so we can generate $2^\gm$ subalgebras by closing
representative elements from  each class.
\end{proof}

\section{ Pr($\gk$) }

\begin{definition} Let $\gk$ be a cardinal.  $Pr(\gk)$ is the following property
 of $\gk$: for all algebras ${\cal M} = (M, {\cal F})$ with $\card{M} = \gk$ and
 $\card{\cal F} < \gk$ there exists $A \subseteq M$ free for $\cal M$ such that
 $\card{A} = \gk$.
\end{definition}

  In some contexts we might wish to make a more complex definition of the
 form ``$Pr(\gk, \gm, D, \gs)$ iff for every algebra on $\gk$ with at most
 $\gm$ functions there is a free set in the $\gs$-complete filter $D$''
 but $Pr(\gk)$ is sufficient for the arguments here.

  We collect some information about $Pr(\gk)$.

\begin{fact} If $\gk$ is Ramsey then $Pr(\gk)$.
\end{fact}

\begin{proof} Let ${\cal M} = (M, {\cal F})$ be an algebra with $\card M = \gk$,
 $\card{\cal F} < \gk$. We lose nothing by assuming $M = \gk$. We can regard $\cal M$
 as a structure for a first-order language $\cal L$ with $\card{\cal L} < \gk$.
 By a standard application of Ramseyness we can get $A \subseteq \gk$ of
 order type $\gk$ with $A$ a set of order indiscernibles for $\cal M$.

   Now $A$ must be free by an easy application of indiscernibility.

\end{proof}

\begin{fact}
\label{cc-fact}
 If $Pr(\gk)$ holds and $\FP$ is a $\gm^+$-c.c.~forcing for
 some  $\gm < \gk$ then $Pr(\gk)$ holds in $V^\FP$.
\end{fact}

\begin{proof} Let $\dot {\cal M} = (\hat \gk, \seq{\dot F_\ga: \ga < \gl})$ name
 an algebra on $\gk$ with $\gl$ functions, $\gl < \gk$. For each $\vec a$ from
 $\gk$ and $\ga < \gl$ we may (by the chain condition of $\FP$)
 enumerate the possible values of $\dot F_\ga(\vec a)$ as
 $\seq{ G_{\ga, \gb}(\vec a) : \gb < \gm}$.

  Now define in $V$ an algebra ${\cal M}^* = (\gk, \seq{G_{\ga, \gb}: \ga < \gl, \gb < \gm})$
 and use $Pr(\gk)$ to get $I \subseteq \gk$ free for ${\cal M}^*$. Clearly
 it is forced that $I$ be free for $\dot {\cal M}$ in $V^\FP$.
\end{proof} 

\section{Building irredundant subsets}

  In this section we will define a combinatorial principle $(*)$ and show that
 if $(*)$ holds then we get the desired conclusion about Boolean algebras.
 We will also consider a limitation  on the possibilities for
 generalising the result.

\begin{definition}[The principle $(*)$]
  The principle $(*)$ is the conjunction of the following two
 statements:
\begin{enumerate} 
\item[S1.] For all infinite cardinals $\gk$, $2^\gk$ is weakly inaccessible
 and
\[
   \gk \le \gl < 2^\gk \implies 2^\gl = 2^\gk.
\]
\item[S2.] For all infinite cardinals $\gk$, $Pr(2^\gk)$ holds.
\end{enumerate}
\end{definition}

\begin{theorem}
\label{startoBA}
 $(*)$ implies that every infinite Boolean algebra $B$
 has an irredundant subset $A$ such that $2^\card{A} = 2^\card{B}$.
\end{theorem}

\begin{proof} 
   Observe that a Boolean algebra can be regarded as an algebra with
  $\aleph_0$ functions (just take the Boolean operations and
 close under composition).

  Define a closed unbounded class of infinite cardinals by
\[
   C = \setof{\gm}{\hbox{$\exists \gth \; 2^\gth = \gm$ or $\gm$ is strong limit}}
\]
 There is a unique $\gm \in C$ such that $\gm \le \card{B} < 2^\gm$. By the
assumption S1,  $2^\gm = 2^{\card B}$, so it will suffice to find a
free subset of size $\gm$. Since $\gm$ is infinite we can find $B_0$
 a subalgebra of $B$ with $\card{B_0} = \gm$, and it suffices to find
 a free subset of size $\gm$ in $B_0$.

\noindent{Case A:} $\gm = 2^\gth$. In this case $\gm > \go$, and S2 implies
 that $Pr(\gm)$ holds, so that $B_0$ has a free (irredundant) subset
 of cardinality $\gm$.

\noindent{Case B:} $\gm$ is strong limit. In this case Fact \ref{BAfact} implies
 that there is an irredundant subset of cardinality $\gm$.
   
\end{proof}

   At this point we could ask about the situation for more general algebras.
 Could every algebra with countably many functions have a large free subset,
 or a large number of subalgebras?
 There is a result of Shelah (see Chapter III of \cite{Sh-e}) which sheds
 some light on this question.

\begin{theorem} If $\gk$ is inaccessible and not Mahlo (for example if
 $\gk$ is the first inaccessible cardinal) then $\gk \nrightarrow [\gk]^2_\gk$.
\end{theorem}

  From this it follows immediately that if $\gk$ is such a cardinal, then
 we can define an algebra $\cal M$ on $\gk$ with countably many functions
 such that $\card{Sub({\cal M})} = \gk$.  
 
   In the next section we will prove that in the absence of inaccessibles
 some information about general algebras can be extracted from $(*)$. 
 The model of $(*)$ which we eventually construct will in fact contain
 some quite large cardinals; it will be a set model of ZFC in which there
 is a proper class of cardinals $\ga$ which are $\beth_3(\ga)$-supercompact.

%
%
%
%
%
%

\section{General algebras}

  In this section we prove that in the absence of inaccessible cardinals
 $(*)$ implies that certain algebras have a large number of subalgebras.
 By the remarks at the end of the last section, it will follow that
 the restriction that
 no inaccessibles should exist is essential.

\begin{theorem} Suppose that $(*)$ holds and there is no inaccessible
 cardinal. Let $\gk$ and $\gl$ be
 infinite cardinals.  If ${\cal M} = (M, {\cal F})$
 is an algebra such that $\card{M} = \gl$ and $\card{\cal F} = \gk$, and
 $2^\gk \le \gl$,  then $Sub(\cal M)$ has cardinality $2^\gl$.
\end{theorem}

\begin{proof}

 Let $\cal M$, $\gk$ and $\gl$ 
 be as above. Define a class of infinite cardinals
$C$ by
\[
    C = \setof{\gm}{\hbox{$\exists \gth \; 2^\gth = \gm$ or $\gm$ is singular
 strong limit.}}
\]
  There are no inaccessible cardinals, so $C$ is closed and unbounded
 in the class of ordinals.

  Let $\gm$ be the unique element of $C$ such that $\gm \le \gl < 2^\gm$,
 where we know that $\gm$ exists because $\gl \ge 2^\gk$ and $C$ is a
 closed unbounded class of ordinals.
 Since $2^\gk \le \gl < 2^\gm$, $\gk < \gm$. 

  Now let $N$ be a subalgebra of $\cal M$ of cardinality $\gm$; such a
 subalgebra can be obtained by taking any subset of size $\gm$ and
 closing it to get a subalgebra. S1 implies that $2^\gl = 2^\gm$,
 so it suffices to show that the algebra ${\cal N} = (N, {\cal F} \restriction N)$
 has $2^\gm$ subalgebras.

 We distinguish two cases:

\noindent {Case A:} $\gm$ is singular strong limit. In this case we may
 apply Fact \ref{singularfact} to conclude that $\cal N$ has $2^\gm$
 subalgebras.

\noindent{Case B:} $\gm = 2^\gth$. It follows from S2  that $Pr(\gm)$ is
 true, hence there is a free subset of size $\gm$ for $\cal N$.
 But now by Fact \ref{freefact} we may generate $2^\gm$ subalgebras.

 This concludes the proof.

\end{proof}
 
   For algebras with countably many functions the last result guarantees
 that there are many subalgebras as long as $\card{M} \ge 2^{\aleph_0}$.
 In fact we can do slightly better here, using another result of Shelah
 (see \cite{avRuSh}).

\begin{theorem} If $C$ is a closed unbounded subset of $[\aleph_2]^{\aleph_0}$
 then $\card{C} \ge 2^{\aleph_0}$.
\end{theorem}

\begin{corollary} If $(*)$ holds and there are no inaccessibles, then for
 every algebra $\cal M$ with $\card{M} > \aleph_1$ and countably many functions
 $\card{Sub({\cal M})} = 2^\card{M}$.
\end{corollary}

\begin{proof} If $\card{M} \ge 2^{\aleph_0}$ we already have it. If
 $\aleph_2 \le \card{M} < 2^{\aleph_0}$ then apply the theorem we
 just quoted and the fact that $2^\card{M} = 2^{\aleph_0}$.
\end{proof}

\section{How to make a model of $(*)$.}

\label{starsection}

  In this section we sketch the argument of \cite{FW} and
 indicate how we modify it to get a model in which $(*)$
 holds. We begin with a brief review of our conventions
 about forcing. For us $p \le q$ means that $p$ is stronger than
 $q$, a forcing is $\gl$-closed if every decreasing sequence of
 length less than $\gl$ has a lower bound, a forcing is $\gl$-dense
 if it adds no $<\gl$-sequence of ordinals, and $Add(\gk, \gl)$ 
 is the forcing for adding $\gl$ Cohen subsets of $\kappa$.

 Foreman and Woodin begin the construction in \cite{FW}
 with a model $V$ in which $\gk$ 
is supercompact, and in which for each finite $n$ they have arranged that $\beth_n(\gk)$
 is weakly inaccessible and $\beth_n(\gk)^{<\beth_n(\gk)} = \beth_n(\gk)$.
They force with a rather complex forcing $\FP$, and pass to a submodel
 of $V^\FP$ which is of the form $V^{\FP^\gp}$ for $\FP^\pi$ a
 projection of $\FP$. 

  $\FP$ here is a kind of hybrid of Magidor's forcing from 
 \cite{SCH1} to violate the Singular Cardinals Hypothesis at
 $\aleph_\omega$
 and Radin's forcing from \cite{Radin}. Just as in \cite{SCH1} 
 the forcing $\FP$ does too much damage to $V$ and the desired
 model is an inner model of $V^\FP$, but in the context of \cite{FW} it
 is necessary to be more explicit about the forcing for which the
 inner model is a generic extension.

 The following are the key properties of $\FP^\pi$.

\begin{enumerate}
\item[P$1$.] $\gk$ is still inaccessible (and in fact is $\beth_3(\gk)$-supercompact)
 in $V^{\FP^\pi}$. 
\item[P$2$.] $\FP^\pi$ adds (among other things) a generic club of order type $\gk$ in $\gk$.
In what follows we will assume that a generic $G$ for $\FP^\pi$ is
 given, and enumerate this club in increasing order as $\seq{\gk_\ga:\ga<\gk}$.
\item[P$3$.] {\bf In $\bf V$\/} each cardinal $\gk_\ga$ reflects to some extent the properties
 of $\gk$. In particular in $V$ each
 $\gk_\ga$ is measurable, and for each $n$
 we have that $\beth_n(\gk_\ga)$ is weakly inaccessible and $\beth_n(\gk_\ga)^{<\beth_n(\gk_\ga)}
 = \beth_n(\gk_\ga)$. 
\item[P$4$.] For each $\ga$, if $p \in \FP^\pi$ is a condition which determines
 $\gk_\ga$ and $\gk_{\ga+1}$ then $\FP^\pi \restriction p$ (the suborder
 of conditions which refine $p$) factors as
\[
    \FP_\ga \times Add(\beth_4(\gk_\ga), \gk_{\ga+1}) \times \FQ_\ga,
\]
where
\begin{enumerate}
\item $\FP_\ga$ is $\gk_\ga^+$-c.c.~if $\ga$ is limit and
 $\beth_4(\gk_\gb)^+$-c.c.~if $\ga=\gb+1$. 
\item All bounded subsets of $\beth_4(\gk_{\ga+1})$ which occur in
 $V^{\FP^\pi \restriction p}$ have already appeared in the extension by 
$\FP_\ga \times Add(\beth_4(\gk_\ga), \gk_{\ga+1})$.
\end{enumerate}

\end{enumerate}

Given all this it is routine to check that if we truncate $V^{\FP^\pi}$
at $\gk$ then we obtain a set model of ZFC in which GCH fails everywhere,
 and moreover exponentiation follows the pattern of clause 1 in 
$(*)$.

   The reader should think of $\FP^\pi$ as shooting a club of
 cardinals through $\gk$, and doing a certain amount of work
 between each successive pair of cardinals. Things have been
 arranged so that GCH will fail for a long region past each
 cardinal $\gk_\ga$ on the club, so that what needs to be done is to blow up
 the powerset of a well-chosen point in that region to have
 size $\gk_{\ga+1}$.

  For our purposes we need to change the construction of \cite{FW}
 slightly, by changing the forcing which is done between the points
 of the generic club to blow up powersets. The reason for the change
 is that we are trying to get $Pr(\gl)$ to hold whenever $\gl$ is
 of the form $2^\gee$, and in general it will not be possible to
 arrange that $Pr(\gl)$ is preserved by forcing with $Add(\gl, \gr)$.
 
  The solution to this dilemma is to replace the Cohen forcing
 defined in $V$ by a Cohen forcing defined in a well chosen inner
 model. We need to choose this inner model to be small enough that
 Cohen forcing from that model preserves $Pr(\gl)$, yet large enough that
 it retains some degree of closure sufficient to make the arguments
 of \cite{FW} go through. For technical reasons we will be adding subsets
 to $\beth_5(\gk_\ga)$, rather than $\beth_4(\gk_\ga)$ as in the case
 of \cite{FW}.

  To be more precise we will want to replace clauses P$3$ and P$4$ above by
\begin{enumerate}
\item[P$3^*$.] For each $n < 6$ it is the case in $V$ that
 $\gk_\ga$ is measurable, $\beth_n(\gk_\ga)$ is weakly inaccessible, and
 $\beth_n(\gk_\ga)^{<\beth_n(\gk_\ga)} = \beth_n(\gk_\ga)$. 
Also $Pr(\beth_n(\gk_\ga))$ holds for each such $n$.

              There exists for all possible values of $\gk_\ga$, $\gk_{\ga+1}$
 a  notion of forcing $Add^*(\beth_5(\gk_\ga),\gk_{\ga+1})$ (which we will
 denote by $\FR_\ga$ in what follows) such that
\begin{enumerate} 
\item $\FR_\ga$ is $\beth_4(\gk_\ga)$-closed,  $\beth_4(\gk_\ga)^+$-c.c.~and
  $\beth_5(\gk_\ga)$-dense.
\item $\FR_\ga$ adds $\gk_{\ga+1}$ subsets to  $\beth_5(\gk_\ga)$.
\item In $V^{\FR_\ga}$ the property $Pr(\beth_n(\gk_\ga))$ still holds
\footnote{The hard work comes in the case $n=5$, density guarantees that
the property survives for $n<5$.}
 for $n < 6$.
\end{enumerate}
\item[P$4^*.$] If $p$ determines $\gk_\ga$ and $\gk_{\ga+1}$ then
 $\FP^\pi \restriction p$ factors as
\[
    \FP_\ga \times \FR_\ga \times \FQ_\ga,
\]
where
\begin{enumerate}
\item $\FP_\ga$ is $\gk_\ga^+$-c.c. if $\ga$ is limit, and
 $\beth_5(\gk_\gb)^+$-c.c.~if $\ga = \gb+1$.
\item $\FR_\ga =  Add^*(\beth_5(\gk_\ga),\gk_{\ga+1})$ is as above.
\item All bounded subsets of $\beth_5(\gk_{\ga+1})$ in $V^{\FP^\pi \restriction p}$ are
 already in the extension by $\FP_\ga \times \FR_\ga$.
\end{enumerate}
\end{enumerate}

 We'll show that given a model in which P$1$, P$2$, P$3^*$ and P$4^*$ hold we can construct
 a model of $(*)$.

\begin{lemma} Let $G$ be generic for a modified version of $\FP^\pi$ obeying
 P$1$, P$2$, P$3^*$ and  P$4^*$.  If we define $V_1= V[G]$,
then $V_1$ is a model of ZFC in which $(*)$ holds.
\end{lemma}

\begin{proof} S1 is proved just as in \cite{FW}.

   It follows from our assumptions that in $V_1$ we have the following
 situation:
\begin{itemize}
\item $\gm$ is strong limit if and only if $\gm = \gk_\gl$ for some
 limit $\gl$.
\item Cardinal arithmetic follows the pattern that for $n < 5$
\[
    \beth_n(\gk_\ga) \le \gth < \beth_{n+1}(\gk_\ga) \implies 2^\gth = \beth_{n+1}(\gk_\ga),
\]
 while
\[
    \beth_5(\gk_\ga) \le \gth < \gk_{\ga+1} \implies 2^\gth = \gk_{\ga+1}.
\]
\end{itemize}

 We need to show that for all infinite $\gl$ we have $Pr(2^\gl)$.
The proof divides into two cases.

\begin{enumerate}
\item[Case 1:] $2^\gl$ is of the form $\beth_m(\gk_\gee)$ where $1 \le m \le 5$
 and $\gee$ is a limit ordinal. By the factorisation properties  in P$4^*$
 above it will suffice to show that $Pr(2^\gl)$ holds in the extension by
 $\FP_\gee \times \FR_\gee$. 

   Certainly $Pr(2^\gl)$ holds in the extension by $\FR_\gee$, because P$3^*$
  says just that. Also we have that $\gk_\gee^+ < \beth_1(\gk_\gee) \le 2^\gl$,
 and $\FP_\gee$ is $\gk_\gee^+$-c.c.~so that applying Fact \ref{cc-fact}
 we get that $Pr(2^\gl)$ holds in the extension by $\FP_\gee \times \FR_\gee$.

\item[Case 2:] $2^\gl$ is of  the form $\beth_m(\gk_{\gb+1})$ where $0 \le m \le 5$.
Again it suffices to show that $Pr(2^\gl)$ holds in the extension by 
 $\FP_{\gb+1} \times \FR_{\gb + 1}$.

   $Pr(2^\gl)$ holds in $V^{\FR_{\gb+1}}$, $\beth_5(\gk_\gb)^+ < \gk_{\gb+1} \le 2^\gl$
 and $\FP_{\gb+1}$ is $\beth_5(\gk_\gb)^+$-c.c.~so that as in the last case we
 can apply Fact \ref{cc-fact} to get that $Pr(2^\gl)$ in the extension by
  $\FP_{\gb+1} \times \FR_{\gb + 1}$.
\end{enumerate}
\end{proof}

\section{The preparation forcing}

   Most of the hard work in this paper comes in preparing
 the model over which we intend (ultimately) to do our
 version of the construction from \cite{FW}.

  We'll make  use of Laver's ``indestructibility'' 
 theorem from \cite{L} as a labour-saving device. At a certain
 point below we will sketch a proof, since we need a little
 more information than is contained in the statement.

\begin{fact}[Laver]
\label{laverology}
   Let $\gk$ be supercompact, let $\gee < \gk$. Then there
 is a $\gk$-c.c.~and $\gee$-directed closed
 forcing $\FP$,  of cardinality $\gk$, such
 that in the extension by $\FP$ the supercompactness of
 $\gk$ is indestructible under $\gk$-directed closed forcing.
\end{fact}

  Now we will describe how to prepare the model over which we will do
 Radin forcing. Let us start with six supercompact cardinals enumerated
 in increasing order as $\seq{\gk_i: i < 6}$ in some initial model $V_0$. 

\subsection{Step One} We'll force to make each of the $\gk_i$ indestructibly
 supercompact under $\gk_i$-directed closed forcing. To do this let
 $\FS_0 \in V_0$ be Laver's forcing to make $\gk_0$ indestructible under
 $\gk_0$-directed closed forcing, and force with $\FS_0$. In $V_0^{\FS_0}$
  the cardinal $\gk_0$ is supercompact by construction, and the rest  of the $\gk_i$
 are still supercompact because $\card{S_0} = \gk_0$ and supercompactness survives
 small forcing.
  
 Now let $\FS_1 \in V_0^{\FS_0}$ be (as guaranteed by Fact \ref{laverology}) 
 a forcing which makes $\gk_1$ indestructible and is $\gk_0$-directed closed.
 In $V_0^{\FS_0 * \FS_1}$ we claim that all the $\gk_i$ are supercompact and
 that both $\gk_0$ and $\gk_1$ are indestructible. The point for the indestructibility
 of $\gk_0$ is that if $\FQ \in V_0^{\FS_0 * \FS_1}$ is $\gk_0$-directed closed
 then $\FS_1 * \FQ$ is $\gk_0$-directed closed in $V_0^{\FS_0}$, so that $\gk_0$
 has been immunised against its ill effects and is supercompact in 
 $V_0^{\FS_0 * \FS_1 * \FQ}$.

  Repeating in the obvious way, we make each of the $\gk_i$ indestructible.
 Let the resulting model be called $V$. 

\subsection{Step Two} Now we will force over $V$ with a rather complicated
 $\gk_0$-directed closed forcing.
 Broadly speaking we will make $\gk_0$ have the properties that are
 demanded of $\gk_\ga$ in P$3^*$ from the last section.

 Of course $\gk_0$ will still be supercompact
 after this, and we will use a reflection argument to show that we have many
 points in $\gk_0$ which are good candidates to become points on the
 generic club.

 To be more precise we are aiming to make a model $W$ in which the following
 list of properties holds:
\begin{enumerate}
\item[K1.] $\gk_0$ is supercompact.
\item[K2.] For each $i < 5$, $2^{\gk_i} = \gk_{i+1}$ and $\gk_{i+1}$ is
 a weakly inaccessible cardinal with $\gk_{i+1}^{<\gk_{i+1}} = \gk_{i+1}$.
\item[] For any  $\gm > \gk_5$ there is a forcing $\FQ_\gm \subseteq V_\gm$ such
 that (if we define $\FQ_\gee = \FQ_\gm \cap V_\gee$ for $\gee$ 
 with $\gk_5 < \gee < \gm$)
\item[K3.] $\card{\FQ_\gee} = \gee^{<\gk_5}$, $\FQ_\gee$ is $\gk_4$-directed closed,
 $\gk_5$-dense and $\gk_4^+$-c.c.
\item[K4.] In $W^{\FQ_\gee}$,  $2^{\gk_5} \ge \gee$ and $Pr(\gk_i)$
 holds for $i < 6$.
\item[K5.] If $\gee < \gz \le \gm$ then $\FQ_\gee$ is a complete suborder of $\FQ_\gz$.
\end{enumerate}

  For each $i < 5$ let $\FP_i$ be the Cohen conditions (as computed
 by $V$) for adding $\gk_{i+1}$ subsets of $\gk_i$. We will abbreviate
 this by $\FP_i = Add(\gk_i, \gk_{i+1})_V$. Let 
$\FP_{i, j} = \prod_{i \le n \le j} \FP_n$.

  Let $V_1 = V[\dot \FP_{0,3}]$. In $V_1$ we know that $2^{\gk_i} = \gk_{i+1}$
 for $i < 4$, and that $\gk_5$ is still supercompact since $\card{\FP_{0,3}} = \gk_4$.
In $V_1$ define a forcing as in Fact \ref{laverology} for making $\gk_5$ indestructible.
More precisely choose $\FR \in V_1$ such that
\begin{itemize}
\item In $V_1$, $\FR$ is $\gk_5$-c.c.~and $\gk_4^{+50}$-directed closed with cardinality
 $\gk_5$.
\item In $V_1[\dot \FR]$ the supercompactness of $\gk_5$ is indestructible under
 $\gk_5$-directed closed forcing.
\end{itemize}

  Let $V_2 = V_1[\dot \FR]$. Define
$\FQ_\gm = Add(\gk_5, \gm)_{V_2}$, and observe that for every $\gee$ 
 $\FQ_\gee = Add(\gk_5, \gee)_{V_2}$.  This is really the key point in the
 construction, to choose $\FQ_\gm$ in exactly the right inner model (see our remarks
 in the last section just before the definition of P$3^*$).

 Finally let $W = V_2[\dot \FP_4]$, where we stress that $\FP_4$ is Cohen
 forcing as computed in $V$.

We'll prove that we have the required list of facts in $W$.

\begin{enumerate}
\item[K1.] $\gk_0$ is supercompact.
\item[] \begin{proof} 
 Clearly $\FP_{0,3} * \FR$ is $\gk_0$-directed closed in $V$. $\FP_4$ is $\gk_4$-directed
 closed in $V$ so is still $\gk_0$-directed closed in $V_2$, hence
 $\FP_{0,3} * \FR * \FP_4$ is $\gk_0$-directed closed in $V$. By the
 indestructibility of $\gk_0$ in $V$, $\gk_0$ is supercompact in $W$.
\end{proof}
\item[K2.] For each $i < 5$, $2^{\gk_i} = \gk_{i+1}$ and $\gk_{i+1}$ is
 a weakly inaccessible cardinal with $\gk_{i+1}^{<\gk_{i+1}} = \gk_{i+1}$.
\item[] \begin{proof}
 In $V_1$ we have this for $i < 4$. $\FP_4$ is $\gk_4$-closed and $\gk_4^+$-c.c.~in $V$,
 $\FP_{0,3}$ is $\gk_3^+$-c.c.~even in $V[\dot\FP_4]$ so by the usual arguments
 with Easton's lemma
  $\FP_4$ is $\gk_4$-dense and $\gk_4^+$-c.c.~in $V_1$. $\FR$ is $\gk_4^{+50}$-closed
 in $V_1$ so that $\FP_4$ is $\gk_4$-dense and $\gk_4^+$-c.c.~in $V_2$.
 Moreover, by the closure of $\FR$ we still have the claim for $i < 4$ in
 $V_2$. $\gk_5$ is inaccessible (indeed supercompact) in $V_2$ so
 that after forcing with $\FP_4$ we have the claim for $i < 5$.
 \end{proof}
\item[] Recall that $\FQ_\gee = Add(\gk_5, \gee)_{V_2}$.
\item[K3.] $\card{\FQ_\gee} = \gee^{<\gk_5}$, $\FQ_\gee$ is $\gk_4$-directed closed,
 $\gk_5$-dense and $\gk_4^+$-c.c.
\item[] \begin{proof} The cardinality statement is clear. 
 $\FQ_\gee$ is $\gk_5$-directed closed in $V_2$, so that (by what we 
 proved about the properties of $\FP_4$ in $V_2$ during the proof
 of K2) $\FQ_\gee$ is
 $\gk_4$-directed closed and $\gk_5$-dense in $W$. Since $\FQ_\gee$
 is $\gk_5$-closed in $V_2$, $\FP_4$ is $\gk_4^+$-c.c.~in
 $V_2[\dot\FQ_\gee]$. $\FQ_\gee$ is $\gk_5^+$-c.c.~in $V_2$
 so that $\FQ_\gee \times \FP_4$ is $\gk_5^+$-c.c.~in $V_2$,
 so that finally $\FQ_\gee$ is $\gk_5^+$-c.c.~in $W$.
\end{proof}  
\item[K4.] In $W^{\FQ_\gee}$,  $2^{\gk_5} \ge \gee$ and $Pr(\gk_i)$
 holds for $i < 6$.
\item[] \begin{proof} Easily $2^{\gk_5} \ge \gee$. We break up the rest
 of the proof into a series of claims.
\begin{claim} For $i < 4$, $Pr(\gk_i)$ holds in $W[\dot \FQ_\gee]$.
\end{claim}
\begin{proof} $\gk_i$ is supercompact in $V[\dot \FP_{i,3}]$, because 
 $\FP_{i,3}$ is $\gk_i$-directed closed. $\FP_{0,i-1}$ has
 $\gk_{i-1}^+$-c.c.~in that model so that $Pr(\gk_i)$ holds in
 $V_1$. Also we know that $W[\dot \FQ_\gee] = V_1[\dot\FR][\dot\FP_4][\dot\FQ_\gee]$
 is an extension of $V_1$ by $\gk_4$-dense forcing, so $Pr(\gk_i)$ 
 is still true in $W[\dot \FQ_\gee]$.
\end{proof}
\begin{claim} $Pr(\gk_4)$ holds in $W[\dot \FQ_\gee]$. 
\end{claim}
\begin{proof} $W = V_2[\dot\FP_4] = V[\dot \FP_{0,4}][\dot\FR]$.
 Arguing as in the previous claim, $Pr(\gk_4)$ holds in $V[\FP_{0,4}]$.
$\FP_4$ is $\gk_4^+$-c.c. in $V_1$ so $\FR$ is $\gk_4^+$-dense in
 $V[\dot \FP_{0,4}]$,hence $Pr(\gk_4)$ holds in $W$. Finally $\FQ_\gee$
 is $\gk_5$-dense in $W$, so $Pr(\gk_4)$ holds in $W[\dot\FQ_\gee]$.
\end{proof}
\begin{claim} $Pr(\gk_5)$ holds in $W[\dot\FQ_\gee]$. 
\end{claim}
\begin{proof} $W[\dot\FQ_\gee] = V_2[\dot\FP_4][\dot\FQ_\gee] = V_2[\dot\FQ_\gee][\dot\FP_4]$.
 As $\FQ_\gee$ is $\gk_5$-directed closed in $V_2$, $\gk_5$ is supercompact
 in $V_2[\FQ_\gee]$. $\FP_4$ is $\gk_4^+$-c.c.~in $V_2[\dot \FQ_\gee]$,
 so that $Pr(\gk_4)$ is still true in $W[\dot\FQ_\gee]$. 
\end{proof}

  This finishes the proof of K4.

\end{proof}
\item[K5.] If $\gee < \gz \le \gm$ then $\FQ_\gee$ is a complete suborder of $\FQ_\gz$.
\item[] \begin{proof} This is immediate by the uniform definition of
 Cohen forcing.
\end{proof}

\end{enumerate}

   We finish this section by proving that in $W$ we can find a highly
 supercompact
 embedding $j:W \lra M$ with critical point $\gk_0$ such that in $M$
 the cardinal
 $\gk_0$ has some properties resembling K1--K5 above. The point
 is that there will be many $\ga < \gk_0$ which have these properties
  in $W$, and eventually we'll ensure that every candidate to be on
 the generic club added by the Radin forcing has those properties.
 
   For the rest of this section we will work in the model $W$ unless
 otherwise specified.

\begin{definition} A pair of cardinals $(\ga, \gk)$ is {\it sweet\/} 
 iff (setting $\ga_n = \beth_n(\ga)$)
\begin{enumerate}
\item $\ga$ is measurable.
\item For each $i < 6$, $\ga_i$ is weakly inaccessible and 
 $\ga_i^{<\ga_i} = \ga_i$.
\item[] There is a forcing $\FQ^\ga_\gk \subseteq V_\gk$ 
 such that (setting $\FQ^\ga_\gee = \FQ^\ga_\gk \cap V_\gee$)
\item $\card{\FQ^\ga_\gee} = \gee^{<\ga_5}$, $\FQ^\ga_\gee$ is
 $\ga_4$-directed closed, $\ga_5$-dense and $\ga_4^+$-c.c.
\item In the generic extension by $\FQ^\ga_\gee$, $2^{\ga_5} \ge \gee$
 and $Pr(\ga_i)$ holds for $i < 6$. 
\item If $\gee < \gz \le \gk$ then $\FQ^\ga_\gee$ is a complete suborder
 of $\FQ^\ga_\gz$.
\end{enumerate}
\end{definition}

\begin{theorem} Let $\gl = \beth_{50}(\gk_5)$. 
In $W$ there is an embedding $j:W \lra M$ such that
 $\crit(j) = \gk_0$, ${}^\gl M \subseteq M$, and in the model
$M$ the pair $(\gk_0, j(\gk_0))$ is sweet.
\end{theorem}

\begin{proof} Recall that we started this section with a model
 $V_0$, and forced with some $\FS_0$ to make $\gk_0$ indestructible.
 Let $\FS$ be the forcing that we do over $V_0^{\FS_0}$ to get
 to the model $V_2$, and recall that we force over $V_2$ with
 $\FP_4$ to get $W$. $\FS * \FP_4$ is $\gk_0$-directed closed
 in $V_0^{\FS_0}$.

  Fix a cardinal $\gm$ much larger than $\gl$, such that $\gm^{<\gk} = \gm$.
  As in \cite{L}, fix in $V_0$ an embedding $j_0:V_0 \lra N$ such that
\begin{enumerate}
\item $\crit(j_0) = \gk_0$, $\gm < j_0(\gk_0)$.
\item $V_0 \models {}^\gm N \subseteq N$.
\item $j_0(\FS_0) = \FS_0 * \FS * \FP_4 * \FR$, where $\FR$ is
 $\gm^+$-closed in $N^{\FS_0 * \FS * \FP_4}$. 
\end{enumerate}
  Let $G_0$ be $\FS_0$-generic over $V_0$, let $g_1 * g_2$ be
 $\FS * \FP_4$-generic over $V_0[G_0]$. Let $W = V_0[G_0][g_1][g_2]$
 and let $N^+ = N[G_0][g_1][g_2]$.

 Now $N^+$
 is closed under $\gm$-sequences in $W$ so that
 the factor forcing $\FR$ is $\gm^+$-closed in $W$.
 Let $H$ be $\FR$-generic over  $W$, then it is easy
 to see that $j_0$ lifts to 
\[
      j_1:V_0[G_0] \lra N^+[H].
\]
  Now we know that $j_1 \restriction \FS * \FP_4 \in N[G_0]$,
  and $g_1 * g_2 \in N^+[H]$,
 so that $j[g_1 * g_2] \in  N^+[H]$. By the directed
 closure of $j(\FS * \FP_4)$ in  $N^+[H]$ we can find
 a ``master condition'' $p$ (a lower bound in $j(\FS * \FP_4)$ 
 for  $j[g_1 * g_2]$) and then force over $W[H]$ 
 with $j(\FS * \FP_4) \restriction p$ to get a generic $X$ and
 an embedding
\[
      j_2:W \lra N^* = N^+[H][X].
\]

   Of course this embedding is not in $W$ but
 a certain approximation is. Observe that $\FR * j(\FS * \FP_4)$
 is $\gm^+$-closed in $W$. Now factor $j_2$ through
 the ultrapower by 
\[
     U = \setof{A \subseteq \powerset_\gk \gm}{j[\gm] \in j_2(A)},
\]
\noindent to get a commutative triangle
\commtriangle{W}{M}{N^*}{j}{k}{j_2}

    By the closure $U \in W$, so $j$ is an {\bf internal\/}
 ultrapower map. It witnesses the $\gm$-supercompactness
 of $\gk_0$ in the model  $W$, because (as can easily
 be checked) $U$ is a fine normal measure on $\powerset_\gk \gm$.

  It remains to check that $(\gk_0, j(\gk_0))$ is sweet in $M$. We'll actually
 check that $(\gk_0, j_2(\gk_0))$ is sweet in $N^*$, this suffices because
 $\crit(k) > \gm > \gk_0$.
  
    Clauses 1 and 2 are immediate by the closure of $M$ inside the model
 $W$. Now to witness sweetness let us set
\[
   \FQ^{\gk_0}_{j_2(\gk_0)} = Add(\gk_5, j_2(\gk_0))_{V_0[G_0][g_1]} = Add(\gk_5, j_2(\gk_0))_{N[G_0][g_1]}.
\]
  Clause 3 is true in $N^+$ by the same arguments that we used for
 $W$ above. $N^*$ is an extension of $N^+$ by highly closed forcing so that
 Clause 3 is still true in $N^*$.
 Clause 5 is easy so we are left with Clause 4.

   Let $h$ be $\FQ^{\gk_0}_\gee$-generic over $N^*$ for $\FQ^{\gk_0}_\gee$. Then $h$
 is generic over $N^+$. By the chain condition of $\FQ^{\gk_0}_\gee$
 and the closure of $N^+$ in $W$,
 $W$ and $N^+$ see the same set of antichains
 for $\FQ^{\gk_0}_\gee$. Hence $h$ is $\FQ^{\gk_0}_\gee$-generic over $W$,
 and we showed that $Pr(\gk_i)$ holds in $W[h]$.
 But then by closure again $Pr(\gk_i)$ holds in $N^+[h]$,
 by Easton's lemma $H * X$ is generic over  $N^+[h]$ for highly
 dense forcing so that finally $Pr(\gk_i)$ holds in $N^*[h]$.
\end{proof}  

\section{The final model}
 In this section we will indicate how to modify the construction of \cite{FW},
 so as to force over the model $W$ obtained in the last section and produce a
 model of $(*)$. The modification that we are making is so minor that it
 seems pointless to reproduce the long and complicated definitions of
 $\FP$ and $\FP^\pi$ from \cite{FW}; in this section we just give a sketch
 of what is going on in \cite{FW}, an indication of how the construction
 there is to be modified, and then some hints to the diligent reader as to
 how the proofs in \cite{FW} can be changed to work for our modified
 forcing.

   As we mentioned in Section \ref{starsection}, the construction of \cite{FW}
 involves building a forcing $\FP$ and then defining a projection $\FP^\pi$;
 the idea is that the analysis of $\FP$ provides information which shows that
 $\FP^\pi$ does not too much damage to the cardinal structure of the ground model.
 We start by giving a sketchy account of $\FP^\pi$ as constructed in \cite{FW}.

   Recall the the aim of $\FP^\pi$ is to add  a sequence
$\seq{(\gk_\ga, F_\ga) : \ga < \gk}$ such that:
\begin{enumerate}
\item $\vec \gk$ enumerates a closed unbounded subset of $\gk$.
\item $F_\ga$ is $Add(\beth_4(\gk_\ga), \gk_{\ga+1})$-generic over $V$.
\item In the generic extension, cardinals are preserved and $\gk$ is
 still inaccessible.
\end{enumerate}

   A condition will prescribe finitely many of the $\gk_\ga$, and for each
 $\gk_\ga$ that it prescribes it will give some information about $F_\ga$;
 it will also place some constraint on the possibilities for adding in 
 new values of $\gk_\ga$ and for giving information about the corresponding
 $F_\ga$. The idea goes back to Prikry forcing, where a condition prescribes
 an initial segment of the generic $\go$-sequence and puts constraint on
 the possibilities for adding further points.

  To be more specific $\FP^\pi$ is built from a pair $(\vec w, \vec \F)$
 where $w(0) = \gk$ and $w(\ga)$ is a measure on $V_\gk$ for $\ga > 0$.
 $w(\ga)$ will concentrate on pairs $(\vec u, \vec h)$ that resemble
 $(\vec w \restriction \ga, \vec \F   \restriction \ga)$;
 in particular if we let $\gk_{\vec u} =_{def} u(0)$, $\gk_{\vec u}$ will
 be a cardinal that resembles $\gk$.
 $\F(\ga)$
 will be an ultrafilter on a Boolean algebra $\FQ(\vec w, \ga)$,
 which  consists of functions $f$ such that $\dom(f) \in w(\ga)$
 and $f(\vec u, \vec k) \in RO(Add(\beth_4(\gk_{\vec u}), \gk))$, with
 the operations defined pointwise. 

 Conditions in $\FP^\pi$ will have the form 
\[
\seq{
    (\vec u_0, \vec k_0, \vec A_0, \vec f_0, s_0), \ldots,
    (\vec u_{n-1}, \vec k_{n-1}, \vec A_{n-1}, \vec f_{n-1}, s_{n-1}), 
    (\vec w, \vec \F, \vec B, \vec g)
    }
\]
 where 
\begin{enumerate}
\item $(\vec u_i,\vec k_i) \in V_\gk$, $\gk_{\vec u_0} < \ldots \gk_{\vec u_{n-1}} < \gk$.
\item $\dom (f_i(\ga)) = A_i(\ga) \in u_i(\ga)$, $\dom (g(\ga)) = B(\ga) \in w(\ga)$.
\item $f_i(\ga) \in k_i(\ga)$, $g(\ga) \in \F(\ga)$.
\item $s_i \in Add(\beth_4(\gk_{u_i}), u_{i+1})$ for $i < n-1$, 
 $s_{n-1} \in Add(\beth_4(\gk_{u_{n-1}}), \gk)$. 
\end{enumerate}

   The aim of this condition is to force that the cardinals $\gk_{u_i}$
 are on the generic club, and that the conditions $s_i$ are in the corresponding
 generics. It is also intended to force that if we add in a new cardinal $\ga$ 
 and Add condition $s$ between $\gk_{u_i}$ and $\gk_{u_{i+1}}$ then there is
  a pair $(\vec v, \vec h)$ in some $A_{i+1}(\gb)$ such that 
 $\ga = \gk_{\vec v}$ and $s \le k_{i+1}(\gb)(\vec v, \vec h)$. 
 This motivates the ordering (which we do not spell out in detail here);
 a condition is refined either by refining the ``constraint parts''
 or by adding in new quintuples that obey the current constraints
 and which impose constraints compatible with the current ones.

   The key fact about $\FP^\pi$ is that given a condition $p$, any question about the generic
 extension can be decided by refining the constraint parts of $p$ (this is modelled
 on Prikry's well-known lemma about Prikry forcing). It will follow from this that
 bounded subsets of $\gk$ are derived from initial segments of the generic as in
 P4 b) from Section \ref{starsection}, which in turn will imply that $\FP^\pi$
 has the desired effect of causing the GCH to fail everywhere below $\gk$.
 It remains to be seen that forcing with $\FP^\pi$ preserves the large cardinal
 character of $\gk$; this is done via a master condition argument which
 succeeds because $\vec w$ is quite long (it has a so-called {\it repeat point})
 and the forcing $Add(\beth_4(\gl), \gm)$ is $\beth_4(\gl)$-directed-closed.

   We can now describe how we modify the construction of \cite{FW}.
 To bring our notation more into line with \cite{FW}, let $\gk = \gk_0$,
 let $V=W$, let $j:V \lra M$ be the embedding constructed at the
 end of the last section. For each $\ga < \gk$ such that $(\ga, \gk)$
 is sweet let us fix $\seq{\FQ^\ga_\gee: \gee \le \gk}$ witnessing this.
 Our modification of \cite{FW} is simply to restrict attention to those
 pairs $(\vec u, \vec k)$ such that $(\gk_{\vec u}, \gk)$ is sweet (there
 are sufficiently many because $(\gk, j(\gk))$ is sweet in $M$) and
 then to replace the forcing $Add(\beth_4(\gk_{\vec u}), \gk)$ by
 the forcing $\FQ^{\gk_{\vec u}}_\gk$. It turns out that the proofs
 in \cite{FW} go through essentially unaltered, because they only
 make appeal to rather general properties (chain condition, distributivity,
 directed closure) of Cohen forcing.

 There is one slightly subtle point,
 which is that is we will need a certain uniformity in the dependence of
 $\FQ^\ga_\gee$ on $\gee$, as expressed by the equation
 $\FQ^\ga_\gee = \FQ^\ga_\gz \cap V_\gee$.
 The real point of this is that when we consider a condition in $j(\FP^\pi)$ 
 which puts onto the generic club a point $\ga < \gk$ followed by the point $\gk$,
 the forcing which is happening between $\ga$ and $\gk$ is $j(\FQ)^\ga_\gk$;
 it will be crucial that $j(\FQ)^\ga_\gk = j(\FQ^\ga_\gk) \cap V_\gk = \FQ^\ga_\gk$,
 that is the same forcing that would be used in $\FP^\pi$ between $\ga$ and $\gk$.

  We conclude this section with a short discussion of the necessary changes in the 
 proofs. We assume that the reader has a copy of \cite{FW} to hand; all references
 below are to theorem and section numbers from that paper.

\medskip \noindent  Section Three:  We will use $j:V \lra M$ to construct the master sequence $(\vec M, \vec g)$.
 Recall that we chose $j$ to witness $\beth_{50}(\gk_5)$-supercompactness
 of $\gk_0$, it can be checked that this is enough to make all the 
 arguments below go through.

    When we build $(\vec M, \vec g)$ we will choose each $g_\ga$ so
 that $\dom (g_\ga)$ contains only $(\vec u, \vec h)$ with
 $(\gk_{\vec u}, \gk)$ sweet, and then  let 
  $g_\ga(\vec u, \vec h) \in \vec \FQ^{\gk_{\vec u}}_\gk$. As
 in \cite{FW}, if $j_\ga: V \lra N_\ga$ is the ultrapower by $M_\ga$
 then $W \models {}^{\gk_3} N_\ga \subseteq N_\ga$. If we let $F$ be
 the function given by
\[
     F: (\vec u, \vec h) \longmapsto \FQ^{\gk_{\vec u}}_\gk
\]
 then $[F]_{M_\ga}$ is $\gk_4$-closed in $N_\ga$, hence is $\gk_3^+$-closed
 in $V$.

  This closure will suffice to make appropriate versions of 3.3, 3.4 and
 3.5 go through.

\medskip \noindent  Section Four:  4.1 is just as in \cite{FW}. A version of 4.2 goes through because
 $\FQ^{\gk_a}_\gk \subseteq V_\gk$.

\medskip \noindent  Section Five: In the definition of ``suitable quintuple'' we demand that
 $k_\gd(b) \in \FQ^{\gk_b}_{\gk_u}$ and $s \in \FQ^{\gk_u}_\gk$.

 The definition of ``addability'' is unchanged. Clauses d) and e)
 of that definition still make sense because (e.g.) $\FQ^{\gk_a}_{\gk_v}$
 is a complete suborder of $\FQ^{\gk_a}_{\gk_w}$. 

 Versions of 5.1 and 5.2 go through without change.

  In the definition of ``condition in $\FP$'' we demand that
 $k_\gd(b) \in \FQ^{\gk_b}_{\gk_w}$, $s \in \FQ^{\gk_{u_i}}_{\gk_{u_{i+1}}}$.
 The ordering on $\FP$ is as before. 

  The definitions of ``canonical representatives'' and of ``upper and lower parts''
are unchanged, as is 5.3. 

  We can do a ``local'' version of our argument for 3.3 above to show 5.4
 goes through. Similarly 5.5, 5.6 and 5.8 go through.

\medskip \noindent  Section Six:  We need to make the obvious changes in the definitions of ``suitable quintuple''
 and of the projected forcing $\FP^\pi$. The proofs all go through because the forcing
 $\FQ^\ga_\gee$ has the right chain condition and closure.

\medskip \noindent  Section Seven: The master condition argument of
this section goes through  because the forcing notion $j(\ga
\longmapsto \FQ^\ga_\gk)(\gk)$ is $\gk_4$-directed closed. 

\section{Conclusion}
Putting together the results of the preceding sections we get the following
 result.

\begin{theorem}[Main Theorem] If Con(ZFC + GCH + six supercompact cardinals) then 
 Con(ZFC + every infinite Boolean algebra $B$ has an irredundant subset
 $A$ such that $2^\card{A} = 2^{\card B}$.
\end{theorem}

   Shelah has shown that the conclusion of this consistency result implies
 the failure of weak square, and therefore needs a substantial large cardinal
 hypothesis.

\end{document}